\documentclass[12pt]{article}
\usepackage{amsfonts, amsmath,euscript}

\author{Gleb G. Gusev\thanks{Partially supported by the grants
RFBR-10-01-00678, RFBR-08-01-00110-a, RFBR and SU HSE 09-01-12185-off-m, NOSH-8462.2010.1.}}
\title{Monodromy zeta-function of a polynomial on a complete intersection and Newton polyhedra}
\date{}

\DeclareMathOperator{\pt}{pt}

\DeclareMathOperator{\Vol}{Vol}

\newtheorem{thm}{Theorem}
\newtheorem{prop}{Proposition}
\newtheorem{lem}{Lemma}

\newenvironment{dfntn} {\smallskip\noindent{\bf
Definition\/}.}{\smallskip\par}

\newenvironment{exmpl} {\smallskip\noindent{\bf Example\/}.}{\smallskip\par}
\newenvironment{rmrk} {\smallskip\noindent{\bf Remark\/}.}{\smallskip\par}
\newenvironment{prf} {\noindent{\em Proof\/}.}{{ $\Box$}\smallskip\par}
\newenvironment{prfthm} {\noindent{\em Proof of the theorem\/}.}{{ $\Box$}\smallskip\par}

\renewcommand{\a}{\alpha }
\renewcommand{\d}{\delta }
\newcommand{\D}{\Delta }
\renewcommand{\l}{\lambda }
\renewcommand{\L}{\Lambda }
\renewcommand{\P}{\Pi }
\newcommand{\s}{\sigma }
\newcommand{\p}{\pi}
\newcommand{\OO}{\mathcal{O}}

\newcommand{\xx}{\mathbf x}
\newcommand{\zz}{\mathbf z}
\newcommand{\uu}{\mathbf u}
\newcommand{\kk}{\mathbf k}
\newcommand{\R}{\mathbb{R}}
\newcommand{\C}{\mathbb{C}}
\newcommand{\CP}{\mathbb{CP}}
\newcommand{\DD}{\mathbb{D}}
\newcommand{\ZZ}{\EuScript{Z}}
\newcommand{\Z}{\mathbb{Z}}
\newcommand{\V}{\overline{V}}
\newcommand{\Zeta}{\widetilde{\zeta}}
\newcommand{\Id}{\mathop{\mathrm{Id}}\nolimits}

\begin{document}
\maketitle

\abstract{For a generic (polynomial) one-parameter
deformation of a complete intersection, there is defined its
monodromy zeta-function. We provide explicit formulae for this zeta-function in terms of the corresponding Newton polyhedra in the case the deformation is non-degenerate with respect to its Newton polyhedra. Using this result we obtain the formula for the monodromy zeta-function at the origin of a polynomial on a complete intersection, which is an analog of the Libgober--Sperber theorem.}

\section{Introduction}

Let $\,F_0, F_1,\ldots, F_k\,$ be a set of functions on $\C^n$
defined as polynomials in $n$ complex variables $\,\zz = (z_1,z_2,\ldots, z_n)\,$. Consider the family of varieties
$$
V_c
= \{\zz\in \C^n \mid F_0(\zz)=c,\,F_i(\zz)=0,\, i=1,2,\ldots
,k\},
$$
where $\,c\in \C \,$ is a complex parameter. This family provides the fibration over the punctured neighbourhood of the origin in the parameter space with the fiber $V_c$ over a point $c$ (see below). In this paper we obtain a formula for the monodromy zeta-function of the above fibration in terms of the Newton polyhedra of the polynomials $\,F_0, F_1, \ldots, F_k.\,$ One can consider this result as a global analog of \cite[Theorem 2.2]{G-2} and an analog of {\cite[Theorem 5.5]{Tak}}, where the monodromy zeta-function at infinity is calculated. In section \ref{sectdefor} we consider the case $\,F_0(\zz) = z_n,\,$ where the fibration corresponds to a polynomial deformation of a set of polynomials in $n-1$ variables $\,z_1,z_2,\ldots, z_{n-1}.\,$ The general case is deduced from the partial in section \ref{sectpolyn}. The study is partially motivated by the results of D. Siersma and M. Tibar (\cite{ST}).

Let $\,A = \C^n\setminus Y\,$ be the complement to an arbitrary algebraic hypersurface $Y\subset \C^n$. Let $\,Z =
\{\zz\in \C^n\mid F_i=0,\, i=1,2,\ldots, k\} \cap A$. Denote by $\DD_r$ and $\DD^*_r$ the closed disk in $\C$ of radius $r$ with the centre at the origin and the punctured disk $\,\DD^*_r :=
\DD_r\setminus \{0\}\,$ respectively. It follows from \cite[Theorem 5.1]{VarT} that there exist a finite set $B\subset \C$ such that the restriction
$F = F_0|_Z$ of the function $F_0$ is a topological fibration over $\,\C
\setminus B\,$. In particular, the map $\,F|_{F^{-1}(\DD^*_\d)}\,$ ($ F|_{F^{-1}(\C \setminus \DD_d)}\,$) is a fibration for a small enough $\d$ (for a large enough $d$). Consider the restriction of this fibration to the cycle $\,\{c \cdot \exp(2\p it)\mid t\in [0,1]\},\,$ where $|c|$ is small enough (large enough respectively). Consider the monodromy transformation $\,h_{F,0} \colon Z_c \to Z_c\,$ ($\,h_{F, \infty}
\colon Z_c \to Z_c\,$) of the fiber $Z_c$ over the point $c$ of the result fibration.

The \emph{zeta-function} of an arbitrary transformation $\,h\colon X\to X\,$ of a topological
space $X$ is the rational function:
$$
\zeta_h(t) = \prod_{i\geq 0}(\det(\Id - th_*|_{H^c_i(X;
\C)}))^{(-1)^i},
$$
where $\,H^c_i(X; \C)\,$ denoted the $i$-th homology group with closed support.

\begin{dfntn}
The \emph{monodromy zeta-function (at the origin)} of the function $F_0$ on the set $Z$ is the zeta-function of the transformation $h_{F,0}$, $\zeta_{F_0, Z}(t): = \zeta_{h_{F,0}}(t)$. The \emph{monodromy zeta-function at infinity} of the function $F_0$ on the set $Z$ is the zeta-function of the transformation $h_{F,\infty}$,
$\zeta^{\infty}_{F_0, Z}(t): = \zeta_{h_{F, \infty}}(t)$.
\end{dfntn}

Let $\,S_1, S_2, \ldots, S_n \subset \R^n\,$ be a set of convex bodies. Denote by $\,S_1 S_2\ldots S_n\,$ their Minkovskian mixed volume (see e.g. \cite{Bus}). If $S_j=\emptyset$ for some $j$ we put $\,S_1 S_2\ldots S_n=0.\,$ For a homogenous polynomial $\,T(x_1, x_2, \ldots,
x_k) = \sum \a_{i_1 i_2\ldots i_n}\,x_{i_1} x_{i_2} \ldots
x_{i_n}\,$ of degree $n$, we define $\,T(S_1, S_2,\ldots, S_k\,)$ as $\sum \a_{i_1 i_2\ldots i_n}\,S_{i_1} S_{i_2} \ldots S_{i_n}$.

Let $\,S_1, S_2, \ldots, S_l \subset L \subset \R^n,\,$ be a set of convex bodies that lie in an {$l$-dimensional} rational affine subspace $L$. We define $\,S_1 S_2 \ldots S_l\,$ as the $l$-dimensional \emph{integer mixed volume}, that is, the Minkovskian mixed volume in the affine subspace $L$ normalized in such way that the $l$-dimensional volume of the minimal parallelepiped with integer vertices equals one.

In this paper we obtain a formula for the zeta-function $\,\zeta_{F_0, V}(t),\,$ $V=\{\zz \in \C^n \mid F_1(\zz)=F_2(\zz)=\ldots=F_k(\zz)=0\}\,$ for a generic set of polynomials $\,F_0, F_1,\ldots, F_k\,$ in terms of the integer mixed volumes of the faces of their Newton polyhedra $\,\D_0,\D_1, \ldots, \D_k$.

\section{Zeta-function of a polynomial deformation}\label{sectdefor}
In this section we study the case $\,F_0(\zz) = z_n.\,$ Consider the set of deformations $\,f_{i,\s}(z_1,\ldots, z_{n-1}) := F_i(z_1,\ldots, z_{n-1}, \s)\,$
of the functions $\,f_i := f_{i,0}\,$ on the set $\C^{n-1}$,
$\,i=1,2,\ldots, k,\,$ where $\,\s\in \C\,$ is the parameter of the deformation. The fibre over the point $c$ of the fibration provided by the function $F_0$ on the set $\,\{F_1 = F_2 =\ldots = F_k= 0\}\,$ is $\{f_{1,c}=f_{2,c}=\ldots = f_{k,c}=0\}\times \{c\}$. This fact motivates the following definition.

\begin{dfntn}
 Consider $\,V=\{\zz \in \C^n \mid F_1(\zz)=F_2(\zz)=\ldots=F_k(\zz)=0\}.\,$ We call the zeta-function $\,\zeta_{z_n, V}(t)\,$ $\,(\zeta^{\infty}_{z_n, V}(t))\,$ \emph{the monodromy zeta-function (at infinity) of the deformation} $\,\{f_{i,\s}\mid i=1, 2, \ldots, k\}$.
\end{dfntn}

\subsection{Formulae for the zeta-function of a deformation}\label{sectdeforformulae}
 One has $\,F_i = \sum_{\kk \in \Z^n} F_{i,\kk} \zz^\kk,\,$ where $F_{i,\kk}\in \C,\,\kk\in \Z^n$ are the coefficients of the polynomial $F_i$ and $\,\kk = (k_1, k_2, \ldots, k_n)\,$ are the coordinates in the space $\R^n$ that corresponds to the variables $\,(z_1, z_2,\ldots, z_n)$. Denote by $\,\D_i = \D(F_i)\,$ the \emph{Newton polyhedron} of the polynomial $\,F_i,\,i=1,2, \ldots, k\,$, i.e. the convex hull of the set $\,\{\kk\in \Z^n\mid F_{i,\kk} \neq 0\}.\,$ A subset $I$ of the set $\,\{1, 2,\ldots, n\}\,$ will be called an \emph{index set}. Denote $\,\R^I =\{\kk\in \R^n \mid k_i=0,\,i\notin
I\}.\,$ Let $\,j^I_1<j^I_2<\ldots <j^I_{k(I)}\,$ be the sorted elements of the set $\,\{j\in \{1,2,\ldots,k\}\mid \D_j\cap\R^I \neq \emptyset\}.\,$ Denote $\D_i^I=\D_{j^I_i}\cap \R^I,\,i=1,2,\ldots,k(I)\,$ and $\,F_i^{I} = \sum_{\kk \in \D_{i}^{I}} F_{j^I_i, \kk} \zz^\kk.\,$

An integer covector is called \emph{primitive} if it is not a multiple of another integer covector. Denote by
$\ZZ^I$ the set of primitive covectors in the dual space ${(\R^I)}^*$. For a convex set  $\,S\subset \R^I\,$ and an covector $\,\a \in \ZZ^I$ denote by $S^{\a}$ the subset of $S$ that consists of the points, where the function $\a|_S$ reaches its minimal value: $\,S^{\a} = \{\xx\in S \mid \a(\xx)=\min(\a|_S)\}$. For an arbitrary polynomial $\,P=\sum_{\kk\in \D} P_{\kk} \zz^\kk\,$ with the Newton polyhedron $\,\D \subset \R^I\,$ and a covector $\,\a\in \ZZ^I,\,$ denote by $P^{\a}$ the polynomial $\sum_{\kk\in \D^{\a}} P_{\kk} \zz^{\kk}$. For in index set $I$ that contain $n$, denote by $\,\ZZ^I_+ \subset \ZZ^I\,$ ($\ZZ^I_- \subset \ZZ^I\,$) the subset of covectors $\,\a = \ldots + \a_n\,dk_n,\,$ that have the strictly positive last component: $\,\a_n > 0\,$ (the strictly negative last component: $\,\a_n < 0$).

\begin{dfntn}
Consider a covector $\a\in \ZZ^{\{1,2,\ldots,n\}}$. We say that a system of polynomials  $\,F_1,F_2,\ldots,F_k\,$ is  \emph{$\a$-non-degenerate} with respect to its Newton polyhedra $\,\D_1, \D_2,\ldots, \D_k\,$ if the
$1$-forms $\,dF_i^{\a},\,$ $i=1,2,\ldots, k\,$ are linear independent at all the points of the set $\,\{\zz \in (\C^*)^n \mid F_1^{\a}(\zz)=F_2^{\a}(\zz)=\ldots=F_{k}^{\a}(\zz)=0\}$.

We say that a system of polynomials $\,F_1, F_2,\ldots,F_k\,$ is \emph{$\s$-non-degenerate (at infinity)} with respect to its Newton polyhedra if for each index set $I$ that contain $n$ and for each covector $\,\a \in
\ZZ^I_+\,$ ($\a \in \ZZ^I_-$) the system of polynomials $F^I_1,F^I_2,\ldots, F_{k(I)}^I$ is $\a$-non-degenerate with respect to its Newton polyhedra.

Finally, a system of polynomials $\,F_1, F_2,\ldots, F_k\,$ is called \emph{non-degenerate} with respect to its Newton polyhedra if for each index set $I$ and for each $\,\a \in \ZZ^I\,$ the system of polynomials $F^I_1,F^I_2,\ldots, F_{k(I)}^I$ is $\a$-non-degenerate.
\end{dfntn}

For each index set $\,I\subset \{1,2,\ldots,
n\}\,$ that contain $n$, we define the following rational functions:
\begin{equation*}
\zeta_{\D_1,\D_2,\ldots,\D_k}^I (t)= \prod_{\a\in \ZZ^I_+}
(1-t^{\a(\frac{\partial}{\partial k_n})})^{l!\,Q^l_{k(I)}\left
(\D_{1}^{I,\a}, \D_{2}^{I,\a},\ldots,
\D_{k(I)}^{I,\a}\right )},
\end{equation*}
\begin{equation*}
\zeta_{\D_1,\D_2,\ldots,\D_k}^{I,\infty} (t)= \prod_{\a\in
\ZZ^I_-} (1-t^{-\a(\frac{\partial}{\partial
k_n})})^{l!\,Q^l_{k(I)}\left (\D_{1}^{I,\a}, \D_{2}^{I,\a},\ldots,
\D_{k(I)}^{I,\a}\right )},
\end{equation*}
where $\,l = |I|-1,\,$ $\frac{\partial}{\partial k_n}$ is the vector in $\R^I$ whose only non-zero coordinate is $\,k_n = 1,\,$ and $\,Q^l_k(x_1, x_2,\ldots, x_k): = \left [\prod_{i=1}^{k} \frac{x_i}{1+x_i} \right
]_l,\,$ where $\,[\cdot]_l\,$ denotes the homogenous part of degree $l$ of the power series under consideration. In particular, $\,Q_0^l\equiv 0\,$ for $l>0$ and $Q_0^0 \equiv 1$.

\begin{thm}\label{thm1}
Let a system of polynomials $\,F_1, F_2, \ldots, F_k\,$
be $\s$-non-degenerate with respect to its Newton polyhedra
$\,\D_1, \D_2,\ldots, \D_k.\,$ Then one has
\begin{equation}
\label{1}\zeta_{z_n, V\cap(\C^*)^n}(t) =
\zeta_{\D_1,\D_2,\ldots,\D_k}^{\{1,2,\ldots, n\}} (t),
\end{equation}
\begin{equation}
\label{1a}\zeta_{z_n, V}(t) = \prod_{I\colon n\in I\subset \{1, 2, \ldots, n\}}
\zeta_{\D_1,\D_2,\ldots,\D_k}^I (t),
\end{equation}
where $V=\{\zz \in \C^n \mid
F_1(\zz)=F_2(\zz)=\ldots=F_k(\zz)=0\}\,$.
\end{thm}

\begin{thm}
Let a system of polynomials $\,F_1, F_2, \ldots, F_k\,$
be $\s$-non-degenerate at infinity with respect to its Newton polyhedra $\,\D_1, \D_2,\ldots, \D_k.\,$ Then one has
\begin{equation}
\label{2}\zeta^{\infty}_{z_n, V\cap(\C^*)^n}(t) =
\zeta_{\D_1,\D_2,\ldots,\D_k}^{\{1,2,\ldots, n\}, \infty} (t),
\end{equation}
\begin{equation}
\label{2a}\zeta^{\infty}_{z_n, V}(t) = \prod_{I\colon n\in I\subset
\{1, 2, \ldots, n\}} \zeta_{\D_1,\D_2,\ldots,\D_k}^{I, \infty}
(t),
\end{equation}
where $V=\{\zz \in \C^n \mid
F_1(\zz)=F_2(\zz)=\ldots=F_k(\zz)=0\}$.
\end{thm}

\begin{rmrk}
In the case $k=1$ the equation (\ref{1}) implies:
\begin{equation}
\zeta_{z_n, V\cap(\C^*)^n}(t) = \prod_{\a\in \ZZ^{I_0}_+}
(1-t^{\a(\frac{\partial}{\partial
k_n})})^{(-1)^n (n-1)!\,\Vol_{n-1}(\D_1^{\a})},
\end{equation}
where $\,\Vol_l(\cdot)\,$ denotes the $l$-dimensional integer volume, $\,I_0=\{1,2,\ldots, n\}\,$.
This equality is similar to the equation \cite[Theorem 2.2, (1)]{G-2} for the zeta-function of a singularity deformation.
In fact, denote by $\,f_\s\,$ the germ at the origin of the deformation defined as follows: $\,f_\s(z_1,\ldots, z_{n-1}) = F_1(z_1,\ldots, z_{n-1}, \s).\,$
Using the equation from \cite{G-2} we obtain:
$\,\zeta_{f_\s|_{(\C^*)^{n-1}}}(t) = \prod_{\a\in \ZZ^{I_0}_{++}}
(1-t^{\a(\frac{\partial}{\partial
k_n})})^{(-1)^n (n-1)!\,\Vol_{n-1}(\D_1^{I_0,\a})},\,$ where
$\,\ZZ^{I_0}_{++}\,$ is the subset of covectors in $\ZZ^{I_0}_+$, whose all the components are strictly positive. Hence, the local zeta-function $\,\zeta_{f_\s|_{(\C^*)^{n-1}}}(t)\,$ is a  "natural" factor of the global one $\,\zeta_{z_n,
V\cap(\C^*)^n}(t).\,$ The same observation follows from the localization principle (see on the p.~\pageref{locprinc} and in \cite{GZ-2}).
\end{rmrk}

\begin{exmpl}
Let us assume $\,n=2,\,k=1.\,$ Consider the polynomial $\,F_1(z_1, z_2) = z_1 + z_2 (1+z_1^2).\,$ The equations (\ref{1}), (\ref{1a}), and the corresponding equations from \cite[Theorem 2.2]{G-2} imply that $\,\zeta_{f_\s|_{\C}}(t) = \zeta_{f_\s|_{(\C^*)}}(t) = (1-t),\,$ $\,\zeta_{z_2, V}(t)= \zeta_{z_2, V\cap(\C^*)^2}(t) = (1-t)^2$. One can obtain the same results with the following argumentations. The global fiber is  $\,V_{\s} = \{z_1 \mid F_1(z_1,
\s) = 0\} = \{\frac{-1 + \sqrt{1-4\s^2}}{2\s}\} = \{x_1(\s),
x_2(\s)\},\,$ where $\,x_1(\s)\approx -\s,\,\,x_2(\s)\approx
-\s^{-1}$ when $|\s| \ll 1$. Thus, it consists of two points, one of them is close to the origin and the second one is close to infinity. The monodromy transformation is the identical map of the fiber on itself, therefore $\,\zeta_{z_2, V}(t)= \det((1-t)\Id) = (1-t)^2.\,$ The local fiber $\,\{f_\s(z_1)=0\} =
\{x_1(\s)\}\,$ consists of one point, thus $\zeta_{f_\s|_{\C}}(t)=(1-t)$.
\end{exmpl}

\subsection{Proofs of theorems}

We reduce the calculation of the zeta-function to the integration with respect to the Euler characteristic. (see., ex., \cite{Viro}), using the following {\it localization principle}.

We recall the notion of the zeta-function as applied to a family of sections of a line bundle over a variety. It was introduced by S.M.~Gusein-Zade and D.~Siersma in (\cite{GZ-2}). Let $W$ be a compact complex analytic variety and let $W_1$ be the complement to a compact subvariety of $W$. Let $L$  be a line bundle over $W$ and $q_\s$ be an analytic in $\,\s\in \C_\s\,$ family of sections of the bundle $L$. Let $U$ be the subset of $\,W_1\times \C_\s,\,$ given by the equation $\,q_\s(x)=0.\,$ The restriction to $U$ of the projection $\,W_1\times \C_\s \to \C_\s\,$ is a fibration over the punctured disc $\,\DD^*_r\subset \C_\s\,$ for $\,|r|\ll 1\,$. The {\it zeta-function of a family of sections $q_\s$ restricted to the set $W_1$} is the zeta-function of the monodromy transformation of the above fibration. We denote it $\zeta_{q_\s|_{W_1}}(t)$.

The fibration $L$ is trivial over a neighbourhood of a point $\,x\in W.\,$ Therefore, using a fixed coordinate system one can consider the family of germs at the point $x$ of sections $q_\s$ as a deformation in the parameter $\s$ of a function germ. We denote by $\,\zeta_{q_\s|_{W_1}, x}(t)$ the zeta-function of the germ at the point $\,\s=0\,$ of the above deformation restricted to the set $W_1$ (see, ex., \cite{G-2}).

\begin{thm}[\cite{GZ-2}, "localization principle"]\label{locprinc}
One has:
$$
\zeta_{q_\s|_{W_1}}(t) = \int_W \zeta_{q_\s|_{W_1},x}(t)\,d\chi.
$$
\end{thm}

Using the Newton polyhedra $\D_1, \D_2, \ldots, \D_k$ of the polynomials $F_1, F_2,\ldots, F_k$ one can construct unimodular simplicial partition $\L$ of the dual space $(\R ^{n})^*$, which is fine enough for the system $\{\D_i\}$ in the sense of \cite{Khov1}. Consider the toroidal compactification $X_{\L}$ of the torus $(\C^*)^n$ that corresponds to the partition $\L$. Remind that the standard action of the torus $(\C^*)^n$ on itself can be extended to the action of the torus on the variety $X_{\L}$. The cones $\,\l \in \L\,$ of the partition are in the one-to-one correspondence with the orbits $\,T_\l\subset X_\L\,$ of this action and the orbit $T_\l$ is isomorphic to $(\C^*)^{n-\dim \l}$. Denote by $X_{\L}'$ the complement in $X_{\L}$ to the torus $\,T_{\{0\}}\cong (\C^*)^n.\,$ Let $\,\V\,$ be the closure of the set $V\cap T_{\{0\}} \subset X_{\L}$, denote $V'= \V \cap X_{\L}'$. We prove the following statement.

\begin{lem}
For a fine enough partition $\L$, one has:
\begin{equation}\label{zeta_s1}
\begin{split}
&\zeta_{z_n, V\cap (\C^*)^n}(t) = \int_{V'} {\zeta_{z_n|_{V\cap
(\C^*)^n},\,x}(t)}\,{d\chi},\\
&\zeta_{z_n, V\cap (\C^*)^n}^{\infty}(t) = \int_{V'}
{\zeta_{z_n|_{V\cap (\C^*)^n},\,x}^\infty (t)}\,{d\chi},
\end{split}
\end{equation}
where for a germ at a point $x\in V'$ of a meromorphic function $f$ on the set $\V$ and for an open subset $A\subset \V$ the expression $\zeta_{f|_A, x}(t)$ ($\zeta_{f|_A, x}^{\infty}(t)$) denotes the local zeta-function (at infinity) of the germ at the point $x$ of the function $f$ restricted to the set~$A$.
\end{lem}

\begin{prf}
One can assume the partition $\L$ to be a subdivision of the standard partition $\P$ of the space $(\R^{n})^*$ corresponding to the $n$-dimensional projective space:  $X_{\P} = \CP^n \supset (\C^*)^n$. Let $p \colon X_{\L} \to \CP^n$ be the map of the toric varieties induced by the refinement $\,\L \prec \P\,$. Consider the family of global sections $s_\s,\,\s\in \C,\,$ of the fibration $\OO(1)$ over $\CP^n$ that is defined by the condition $\,s_\s|_{(\C^*)^n} = z_n - \s.\,$ Denote $\p = p\,\circ \,inj$, where $inj \colon \V \hookrightarrow X_{\L}$ is the inclusion map. Denote by $S_\s = \p^*(s_\s)$ the family of sections of the bundle $\p^*(\OO(1))$ that is pull-back of $s_\s$. In a similar manner, consider a family of sections $\,s'_\s,\,\s\in \C\,$ of the fibration $\OO(1)$ that is defined by the condition $\,s'_\s|_{(\C^*)^n}= 1 - \s z_n\,$ and consider its pull-back $S'_\s = \p^*(s'_\s)$.

Due to simple reformulations one can easily see that:
\begin{equation*}
\zeta_{z_n, V\cap (\C^*)^n}(t) = \zeta_{S_\s|_{V\cap (\C^*)^n}}(t), \quad \zeta_{z_n, V\cap (\C^*)^n}^{\infty}(t) =
\zeta_{S'_\s|_{V\cap (\C^*)^n}}(t),
\end{equation*}
\begin{equation*}
\zeta_{z_n|_{V\cap (\C^*)^n},x}(t) = \zeta_{S_\s|_{V\cap (\C^*)^n},x}(t), \quad \zeta_{z_n|_{V\cap (\C^*)^n},x}^{\infty}(t) =
\zeta_{S'_\s|_{V\cap (\C^*)^n},x}(t).
\end{equation*}
Applying theorem $\ref{locprinc}$ to the families $S_c,
S_c'$ we obtain:
$$
\zeta_{S_\s|_{V\cap (\C^*)^n}}(t)  = \int_{\V} {\zeta_{S_\s|_{V\cap (\C^*)^n}, x}(t)}\,{d\chi} = \int_{\V} {\zeta_{z_n|_{V\cap (\C^*)^n}, x}(t)}\,{d\chi},
$$
$$
\zeta_{S'_\s|_{V\cap (\C^*)^n}}(t)= \int_{\V} {\zeta_{S'_\s|_{V\cap (\C^*)^n}, x}(t)}\,{d\chi}= \int_{\V} {\zeta^{\infty}_{z_n|_{V\cap (\C^*)^n}, x}(t)}\,{d\chi}.
$$
Moreover, it is easy to see that  $\,\zeta_{z_n|_{V\cap (\C^*)^n}, x}(t)= \zeta^{\infty}_{z_n|_{V\cap (\C^*)^n}, x}(t)=1\,$ for $\,x\notin V'.\,$ Therefore, using the multiplicative property of the integration we obtain:

$$
\int_{\V} {\zeta_{z_n|_{V\cap (\C^*)^n}, x}(t)}\,{d\chi} = \int_{V'} {\zeta_{z_n|_{V\cap (\C^*)^n},\,x}(t)}\,{d\chi},
$$
$$
\int_{\V} {\zeta^{\infty}_{z_n|_{V\cap (\C^*)^n}, x}(t)}\,{d\chi} = \int_{V'} {\zeta_{z_n|_{V\cap (\C^*)^n},\,x}^\infty (t)}\,{d\chi}.
$$
\end{prf}

Let $\L_+ \subset \L$ ($\L_- \subset \L$) be the subset of cones $\l \in \L$ that is generated by a set of primitive covectors
$\,\a_1,\a_2, \ldots, \a_l\,$ lying in $\ZZ^{\{1,2,\ldots, n\}} \setminus \ZZ^{\{1,2,\ldots, n\}}_-$ (lying in $\ZZ^{\{1,2,\ldots, n\}} \setminus \ZZ^{\{1,2,\ldots, n\}}_+$ respectively). One can assume that $\L$ is fine enough such that $\L_- \cup \L_+ = \L$.

Consider an arbitrary point $\,x_0\in V'.\,$ It is contained in the torus $T_\l$ that corresponds to some $l$-dimensional cone $\l \in \L$, $\,l<n.\,$ This cone lies on the border of an $n$-dimensional cone $\,\l'\in \L.\,$ Denote by $\,\a_1,\a_2, \ldots, \a_l\,$ the primitive integer covectors that generate the cone $\l$. The cone $\l'$ is generated by the covectors $\,\a_1,\a_2,\ldots,\a_l\,$ and some covectors $\a_{l+1},\a_{l+2},\ldots, \a_n\,$ in addition. Consider the coordinate system $\,\uu =(u_1,u_2,\ldots, u_n)\,$ corresponding to the set of covectors $\,(\a_1,\a_2, \ldots, \a_n).\,$ One has: $\,u_i(x_0)=0,\,i\leq l,\,u_i(x_0)\neq 0, i>l.\,$ We express the monomial $z_n$ as a function $F$ in the variables $\uu$:
$$
F(\uu) = b\cdot\,u_1^{\a_1(\partial/ \partial k_n)}  u_2^{\a_2(\partial/ \partial k_n)} \ldots u_l^{\a_l(\partial/ \partial k_n)},
$$
where $\,b(\uu) = \prod_{j=l+1}^n u_j^{\a_j(\partial/ \partial k_n)},\,b(x_0)\in \C^*.\,$
Now we are ready to calculate the values of the integrands $\,\zeta_{z_n|_{V\cap (\C^*)^n},\,x_0}(t)\,$ and $\,\zeta_{z_n|_{V\cap (\C^*)^n},\,x_0}^\infty(t)\,$ in the two following cases.

\begin{enumerate}
\item Assume that $\,\l \in \L_-.\,$ Then the value of the function $F$ at the point $x$ is not zero and therefore $\,\zeta_{z_n|_{V\cap (\C^*)^n},\,x_0}(t) = 1$. Assume that $\,\l \in \L_+\,$ respectively. Then the point $x_0$ is not a pole of the function $F$ and therefore $\zeta_{z_n|_{V\cap (\C^*)^n},\,x_0}^\infty(t)=1$.

\item Assume that $\,\l \in \L\setminus \L_-\,$ ($\l \in \L\setminus \L_+$). The system of polynomials $F_1, F_2, \ldots, F_k$ is $\s$-non-degenerate (at infinity) with respect to its polyhedra $\D_1, \D_2, \ldots, \D_k$. Therefore one has $\,l+k \leq n\,$ and there exists a coordinate system $\,(u_1, \ldots, u_l, w_{l+1}, \ldots, w_n)\,$ in a neighbourhood of the point $x_0$ such that $\,w_i(x_0)=0,\,$ $\,i = l+1, \ldots, n,\,$ and
\begin{equation}\label{F}
F_i =a_i\,u_1^{m_{i,1}} u_2^{m_{i,2}}\ldots u_l^{m_{i,l}} \cdot w_{n-i+1},\,i=1, 2, \ldots, k,
\end{equation}
where $\,m_{i,j}= \min(\a_j|_{\D_i})\,$ and $a_i$ is a germ of function such that $\,a_i(x_0)\neq 0.\,$ Denote $\,V_{x_0} = V\cap (C^*)^n \cap U.\,$ According to (\ref{F}) one has:
$$
V_{x_0} = \{u_i \neq 0,\,i\leq l;\,w_i=0,\,i>n-k\}\subset U.
$$
Hence,
\begin{equation}\label{g}
\zeta_{z_n|_{V\cap (\C^*)^n},\,x_0}(t) =
\zeta_{g|_{\{u_i \neq 0,\,i\leq l\}},\,0} (t)\quad(\zeta_{z_n|_{V\cap (\C^*)^n},\,x_0}^\infty(t)= \zeta_{g|_{\{u_i \neq 0,\,i\leq l\}},\,0}^\infty (t)),
\end{equation}
where $g$ is the germ of function in the variables $\,(u_1,\ldots, u_l,w_{l+1},\ldots,w_{n-k})\,$ that is given by the the following equation:
$$
g = \prod_{j=1}^l u_j^{\a_j(\partial/
\partial k_n)}\cdot b(u_1,\ldots, u_l,w_{l+1},\ldots,w_{n-k},0,\ldots,0).
$$
Using the Varchenko-type formula for meromorphic functions (see \cite{GZ}) we calculate the right-hand side of (\ref{g}). For $\,l=1,\,$ we obtain:
\begin{equation}\label{zetazeta}
\zeta_{z_n|_{V\cap (\C^*)^n},\,x_0}(t) = 1-t^{\a_1(\partial/ \partial k_n)}\quad (\zeta_{z_n|_{V\cap (\C^*)^n},\,x_0}^\infty(t) = 1-t^{-\a_1(\partial/ \partial k_n)}).
\end{equation}
Finally, the both zeta-functions are trivial if $\,l>1.\,$

\end{enumerate}

Now we specify the only case, where the zeta function $\,\zeta_{z_n|_{V\cap (\C^*)^n},\,x_0}(t)\,$ ($\zeta_{z_n|_{V\cap (\C^*)^n},\,x_0}^\infty(t)$) is not trivial. Namely, one should assume that $\,x_0\in T_\l,\,\l \in \L\setminus \L_-$ ($\l \in
\L\setminus \L_+$) and $\dim \l =1$. Denote $\,\a=\a_1.\,$ The set $T_\l\cap V'$ can by defined in the coordinates $(u_2, \ldots, u_{n+1})$ on the torus $\,T_{\l} = \{u_1=0\}\,$ by the system of equations $\,\{Q_1^{\a} = Q_2^{\a}
=\ldots =Q_k^{\a} = 0\,\}$, where
$$
Q_i^{\a} = \sum_{\kk\in \D^{\a}_i} F_{i,\kk}\,u_2^{\a_2(\kk)} u_3^{\a_3(\kk)}\ldots u_n^{\a_n(\kk)}.
$$
Using the main results of \cite{Khov1}, \cite{Khov2} we obtain:
\begin{equation}\label{chichi}
\chi(T_\l\cap V') = (n-1)!\,Q_k^{n-1}(\D(Q_1^{\a}), \D(Q_2^{\a}),\ldots,
\D(Q_k^{\a})),
\end{equation}
where $\D(\cdot)$ denotes the Newton polyhedron of the Laurent polynomial under consideration. The covectors $\,\a_2,\a_3,\ldots, \a_n\,$ define an isomorphism of the integer latices of the hyperplane $\{\a=0\}\subset \R^n_{\kk}$ and the space $\R^{n-1}$, which contain the polyhedra $\D(Q_i^{\a})$. Under this isomorphism, the polyhedra $\D(Q_i^{\a})$ corresponds to a parallel shifts of the polyhedra $\D_i$. Therefore, the corresponding mixed integer volumes coincide and
\begin{equation}\label{R}
Q_k^{n-1}(\D(Q_1^{\a}), \D(Q_2^{\a}),\ldots, \D(Q_k^{\a})) = Q_k^{n-1}(\D^\a_1, \D^\a_2,\ldots, \D^\a_k).
\end{equation}
The equations (\ref{zetazeta}), (\ref{chichi}), (\ref{R}) imply the following answers:
\begin{equation}\label{zeta_s3}
\begin{split}
& \int_{T_{\l}\cap V'} {\zeta_{z_n|_{V\cap
(\C^*)^n},\,x}(t)}\,{d\chi} =  (1-t^{\a(\frac{\partial}{\partial
k_n})})^{(n-1)!\,Q_k^{n-1}(\D^\a_1, \D^\a_2,\ldots, \D^\a_k)},\\
& \int_{T_{\l}\cap V'} {\zeta_{z_n|_{V\cap
(\C^*)^n},\,x}^\infty(t)}\,{d\chi} =
(1-t^{-\a(\frac{\partial}{\partial k_n})})^{(n-1)!\,Q_k^{n-1}(\D^\a_1,
\D^\a_2,\ldots, \D^\a_k)}.
\end{split}
\end{equation}

One can multiply the equations (\ref{zeta_s3}) over all strata $T_\l \subset
X_{\L}'$ of dimension $(n-1)$ corresponding to tori $\l \in \L\setminus \L_-$ ($\l \in \L\setminus \L_+$), apply (\ref{zeta_s1}) and obtain the required equations (\ref{1}) and (\ref{2}). The equations (\ref{1a}) and (\ref{2a}) follow from (\ref{1}) and (\ref{2}) respectively due to the multiplicative property of zeta-functions.

\section{Zeta-function of a polynomial on a complete intersection}\label{sectpolyn}

In this section we obtain the general formula for the zeta-function at the origin of a polynomial $\,F_0=\sum_{\kk \in \Z^n} F_{0,\kk} \zz^\kk\,$ on the set of common zeroes of a set of polynomials $\,F_1,F_2,\ldots,F_k.\,$ We use notations and definitions introduced in the section \ref{sectdefor}. Let $\D_0$ be the Newton polyhedron of the polynomial $F_0$. For an index set $I$, denote $\,\D_0^I = \D_0\cap \R^I,\,F_0^I = \sum_{\kk \in \D_0^I} F_{0,\kk} \zz^\kk.$

For each index set $\,I\subset \{1,2,\ldots,
n\}\,$ consider the following rational function:

\begin{equation}
\label{zeta^I}\Zeta_{\D_0;\D_1,\ldots,\D_k}^I (t):= \prod_{\a\in
\ZZ^I_{\D_0}} (1-t^{m_{\D_0^I}(\a)})^{l!\,\widetilde{Q}^l_{k(I)+1}\left
(\D_0^{I,\a}, \D_1^{I,\a},\ldots, \D_{k(I)}^{I,\a}\right
)},
\end{equation}
where $\,m_{\D_0^I}(\a) = \min(\a|_{\D_0^I})\,$ is the minimal value of the covector $\a$ on the set $\D_0^I$, $\,\ZZ^I_{\D_0}\,$ is the set of covectors $\,\a \in \ZZ^I\,$ such that $\,\min(\a|_{\D_0^I})>0\,$ (for
$\,\D_0^I = \emptyset,\,$ we put $\,\ZZ^I_{\D_0} = \emptyset\,$) and
\begin{equation}\label{tildeQ}
\widetilde{Q}^l_{k+1}(x_0, x_1,\ldots, x_k)= Q^l_k (x_1, x_2,
\ldots, x_k) - Q^l_{k+1}(x_0, x_1,\ldots, x_k)
\end{equation}
is the homogenous polynomial of degree $\,l:= |I|-1.\,$ The following statement follows from Theorem \ref{thm1} and some observations concerning the formula for the Euler characteristic of a non-degenerate complete intersection, which was obtained in \cite{Khov2}.

\begin{thm}
Let the systems of polynomials $\,F_0, F_1, \ldots, F_k\,$ and $\,F_1,F_2 \ldots, F_k\,$ be non-degenerate with respect to its Newton polyhedra
$\,\D_0, \D_1,\ldots, \D_k\,$ and $\,\D_1,\D_2,\ldots, \D_k\,$ respectively. Then one has:
\begin{equation}
\label{3}\zeta_{F_0, V\cap(\C^*)^n}(t) =
\Zeta_{\D_0;\D_1,\ldots,\D_k}^{\{1,\ldots, n\}} (t),
\end{equation}
\begin{equation}
\label{3a}\zeta_{F_0, V}(t) = \prod_{I\subset \{1, \ldots, n\}:\,
I\neq \emptyset} \Zeta_{\D_0;\D_1,\ldots,\D_k}^I (t),
\end{equation}
where $V=\{\zz \in \C^n \mid
F_1(\zz)=F_2(\zz)=\ldots=F_k(\zz)=0\}\,$ is the set of common zeroes of the system $\,F_1,F_2,\ldots, F_n$.
\end{thm}

\begin{rmrk}
Consider the case $k=0$. Using the equations (\ref{3a}) and (\ref{zeta^I}) one can obtain:
\begin{equation*}
\zeta_{F_0, \C^n}(t) = \prod_{I\neq \emptyset} \prod_{\a\in \ZZ^I_{\D_0}}
(1-t^{m_{\D_0}(\a)})^{(-1)^l l!\,\Vol_l(\D_0^{I,\a})}
\end{equation*}
(here we put $\Vol_0(\pt)=1$). This equation is an analog of the Libgober--Sperber theorem (\cite{Lib}) and (in slightly different form) was obtained by ïîëó÷åíà Y.~Matsui and K.~Takeuchi (\cite[section 4]{Tak}).
\end{rmrk}

\begin{prfthm}
Note that the equation (\ref{3a}) follows from the equation (\ref{3}) due to the multiplicative property of zeta-functions. We prove (\ref{3}).

Consider the system of polynomials $G_1,G_2, \ldots, G_{k+1}$ in $n+1$ variables $(\zz, z_{n+1}) = (z_1, z_2, \ldots, z_{n+1})$ those are given by the following equations:
\begin{equation}\label{F->G}
\begin{split}
&G_i(z_1, z_2\ldots, z_{n+1}) = F_i(z_1, z_2\ldots,
z_n),\,i=1,2,\ldots,k,\\
&G_{k+1}(z_1, z_2\ldots, z_{n+1}) = F_0(z_1, z_2\ldots, z_n) - z_{n+1}.
\end{split}
\end{equation}
Consider the set $W = \{(\zz, z_{n+1})\in
\C^{n+1} \mid G_1(\zz)=G_2(\zz)=\ldots =G_{k+1}(\zz) = 0\}$.
The fibrations defined by the maps
$$
V\cap (\C^*)^n \cap
F_0^{-1}(\DD^*_\d)\stackrel{F_0}{\longrightarrow} \DD^*_\d\,
\mbox{ and }\,W\cap (\C^*)^{n+1} \cap \{ 0<|z_{n+1}| \leq \d\}
\stackrel{z_{n+1}}{\longrightarrow}\DD^*_\d,
$$ are obviously isomorphic, therefore one has:
\begin{equation}\label{zeta}
\zeta_{F_0, V\cap (\C^*)^n}(t) = \zeta_{z_{n+1}, W\cap
(\C^*)^{n+1}}(t).
\end{equation}

The space  $\R^n$ with the coordinates $\,(k_1,k_2,\ldots,k_n)\,$ is enclosed in a standard manner in the space $\R^{n+1}$ with the additional coordinate $k_{n+1}$ that corresponds to the variable  $z_{n+1}$. For $\,i\leq k,\,$ the Newton polyhedra of the polynomials $F_i$ and $G_i$ coincide: $\D(G_i) = \D_i$. The Newton polyhedron of the polynomial $G_{k+1}$ is a cone of integer height $1$ over the Newton polyhedron of the polynomial $F_0$, $\,\D(G_{k+1}) =
C\D_0$.

\begin{prop}\label{prop}
For a system of polynomials $\,F_0, F_1,\ldots, F_k\,$ such that both the system itself and the system of polynomials $\,F_1,F_2,\ldots, F_k\,$ are non-degenerate with respect to its Newton polyhedra, the system of polynomials $\,G_1,G_2, \ldots, G_{k+1}\,$ is also non-degenerate with respect to its Newton polyhedra.
\end{prop}
\begin{prf}
Consider an arbitrary subset $I\subset \{1,2,\ldots,
(n+1)\}$ and an arbitrary covector
$\a \in \ZZ^I$. For $\,n+1\notin I,\,$ the conditions of $\a$-non-degeneracy as applied to the system $\{G_i^I\}$ and to the system $\{F_i^I\}$ are obviously equivalent. Assume that $n+1\in I$. Denote $I' = I\setminus \{n+1\},\,\a'
= \a|_{\R^{I'}}$. Extending the notations of the section \ref{sectdeforformulae} to the system of polynomials $\,G_1,G_2,\ldots, G_{k+1},\,$ we get: $\,k(I)=k(I')+1,\,G_i^{I,\a}(\zz, z_{n+1}) = F_i^{I', \a'}(\zz)$ for $i\leq k(I')$. Three cases are possible.
\begin{enumerate}
\item $\a(\frac{\partial}{\partial k_{n+1}})>
\min(\a'|_{\D^{I'}_0})$. Then $(C\D_0\cap \R^I)^{\a} =
\D_0^{I',\a'},\,G_{k(I)}^{I,\a}(\zz, z_{n+1}) = F_0^{I', \a'}(\zz)$.

\item $\a(\frac{\partial}{\partial k_{n+1}})<
\min(\a'|_{\D^{I'}_0})$. Then $(C\D_0\cap \R^I)^{\a} =
\{\frac{\partial}{\partial k_{n+1}}\},\,G_{k(I)}^{I,\a} = -z_{n+1}$.

\item $\a(\frac{\partial}{\partial k_{n+1}})=
\min(\a'|_{\D^{I'}_0})$. Then $(C\D_0\cap \R^I)^{\a} = C(\D_0^{I',\a'})$ is a cone of integer height $1$ over $\D_0^{I',\a'}$, $G_{k(I)}^{I,\a}(\zz, z_{n+1}) = F_0^{I', \a'}(\zz)-z_{n+1}$.
\end{enumerate}
Using $\a'$-non-degeneracy of the systems $\,F_0,F_1,\ldots, F_k\,$ and $\,F_1,F_2,\ldots, F_k\,$ one can easily verify that the 1-forms $dG^{I, \a}_i,\,i=1,2,\ldots,k(I)\,$ are linear independent
at the points of the set $\{(\zz, z_{n+1})\in
(\C^*)^{n+1}\mid G_1(\zz, z_{n+1}) = G_2(\zz, z_{n+1})=\ldots
=G_{k+1}(\zz, z_{n+1})=0\}$.
\end{prf}

It follows from Proposition \ref{prop} that Theorem \ref{thm1} is applicable to the system of polynomials $\,G_1,G_2,\ldots,G_{k+1}$:
\begin{equation}
\zeta_{z_{n+1}, W\cap (\C^*)^{n+1}}(t)= \prod_{\a\in \ZZ^{I_0}_+}
\left (1-t^{\a\left (\frac{\partial}{\partial k_{n+1}}\right
)}\right )^{n!\,Q_{k+1}^n({(C\D_0)}^{\a}, \D_1^{\a},\ldots,
\D_k^{\a})},
\end{equation}
where $I_0=\{1,2,\ldots,n+1\}$. One can easily see that in the above cases 1,2 the exponent
$\,Q_{k+1}^n({(C\D_0)}^{\a}, \D_1^{\a},\ldots,
\D_k^{\a})\,$ equals~$0$. Therefore, one has:
\begin{equation}\label{zetamm}
\zeta_{z_{n+1}, W\cap (\C^*)^{n+1}}(t)= \prod_{\a\in
\ZZ^{I'_0}_{\D_0}}
(1-t^{m_{\D_0}(\a)})^{n!\,Q_{k+1}^n(C(\D_0^{\a}),
\D_1^{\a},\ldots, \D_k^{\a})},
\end{equation}
where $\,I_0'=\{1,2,\ldots,n\}.\,$ Now the equation (\ref{3}) follows from (\ref{zeta}),
(\ref{zetamm}) and the equality
\begin{multline*}
n!\,Q_{k+1}^n((C(\D_0^{\a}), \D_1^{\a},\ldots, \D_k^{\a}) = (n-1)!\,\widetilde{Q}_{k+1}^{n-1}(\D_0^{\a}, \D_1^{\a},\ldots,
\D_k^{\a}),
\end{multline*}
that is a consequence of the following statement.
\end{prfthm}

\begin{prop}
Let $\,\D_0,\D_1,\ldots, \D_k\,$ be a set of integer polyhedra lying in a rational affine hyperplane $\,L\subset \R^{n+1}.\,$ Let $C\D_0$ be the cone over $\D_0$ with the vertex at some point  $\,v\in \R^{n+1}\,$ that is in the integer distance $1$ from the hyperplane $L$. Then the following relation holds:
\begin{equation}\label{Q-tildeQ}
(n+1)!\,Q_{k+1}^{n+1}((C\D_0, \D_1,\ldots, \D_k) = n!\, \widetilde{Q}_{k+1}^n(\D_0,
\D_1,\ldots, \D_k).
\end{equation}
\end{prop}

\begin{prf}
Chose an affine integer coordinate system $\,\kk=(k_1,k_2,\ldots, k_{n+1})\,$ in the space $\R^{n+1}$ in such a manner that $\,L=\{\kk\in \R^n\mid k_{n+1}=0\}\,$ and $\,v = (0,0,\ldots,1).\,$ Chose Laurent polynomials $\,F_0, F_1, \ldots, F_k\,$ in the variables $\zz = \,(z_1,z_2,\ldots, z_n)\,$ with fixed Newton polyhedra $\,\D_0,\D_1,\ldots, \D_k\,$ in such a way that the systems $\,F_0, F_1, \ldots, F_k\,$ and $\,F_1, F_2, \ldots, F_k\,$ are non-degenerate with respect to its Newton polyhedra in the sense of \cite{Khov1}. One can easily show (see Proposition \ref{prop}) that the system of Laurent polynomials $\,G_1, G_2,\ldots, G_{k+1}\,$ in $n+1$ variables defined in terms of the polynomials $\{F_i\}$ by the equations (\ref{F->G}) is also non-degenerate in the sense of \cite{Khov1} with respect to its Newton polyhedra $\,\D_1, \ldots, \D_k, C(\D_0).\,$ Denote:
\begin{equation}
\begin{split}
&V =\{F_0=F_1=\ldots= F_k=0\}\subset (\C^*)^n,\\
&V_1 =\{F_1=F_2=\ldots= F_k=0\}\subset (\C^*)^n,\\
&W=\{G_1=G_2=\ldots= G_{k+1}= 0\}\subset (\C^*)^{n+1}.
\end{split}\notag
\end{equation}

Applying the results of \cite{Khov2} we find the Euler characteristics of the sets $\,V,V_1,W$:
\begin{equation}\label{chiV}
\chi(V) = n!\,Q_{k+1}^n(\D_0, \D_1,\ldots, \D_k),\,\,\chi(V_1) = n!\,Q_k^n(\D_1,\D_2,\ldots, \D_k),
\end{equation}
\begin{equation}\label{chiW}
\chi(W) = (n+1)!\,Q_{k+1}^{n+1}(C\D_0, \D_1,\ldots, \D_k).
\end{equation}
Consider the projection $\,p\colon (\C^*)^{n+1} \to (\C^*)^n\,$ on the coordinate hyperplane with the coordinates $\,(z_1,z_2,\ldots, z_n).\,$ Its restriction $p|_W$ provides an isomorphism between $W$ and $\,V_1\setminus V.\,$ Thus one has the equality
$$
\chi(W) = \chi(V_1)- \chi(V).
$$
Applying this equation and using (\ref{chiV}), (\ref{chiW}), and (\ref{tildeQ}) we get (\ref{Q-tildeQ}).
\end{prf}

\noindent Gusev G.G., Moscow Institute of Physics and Technology,
Independent University of Moscow.
\newline E-mail: gusev@mccme.ru

\end{document}